\documentclass{amsproc}%
\usepackage{graphicx}
\usepackage{amscd}
\usepackage{amsmath}
\usepackage{amsfonts}
\usepackage{hyperref}
\usepackage{amssymb}%
\newtheorem{theorem}{Theorem}[section]
\newtheorem{corollary}[theorem]{Corollary}

\newtheorem{proposition}[theorem]{Proposition}
\theoremstyle{definition}
\newtheorem{definition}[theorem]{Definition}

\theoremstyle{remark}

\numberwithin{equation}{section}
\begin{document}
\title[Lefschetz numbers in control theory]{Applications of Lefschetz numbers in control theory}
\author{Peter Saveliev}
\address{Department of Mathematics, Marshall University }
\email{saveliev@marshall.edu}

\begin{abstract}
The goal of this paper is to develop some applications of the Lefschetz
coincidence theory techniques in control theory. The topics treated are
existence of equilibria and their robustness, controllability and its robustness.

\end{abstract}
\maketitle

\section{Introduction.}

The goal of this paper is to provide examples of what Lefschetz coincidence
theory can contribute to control theory. We discuss existence of equilibria
and their robustness, controllability and its robustness.

We develop applications some techniques, already available in dynamics, in
control theory. A (discrete) dynamical system on a manifold $M$ is simply a
map $f:M\rightarrow M$ and the next state $f(x)$ of the system depends only on
the current one, $x\in M$. The equilibrium set $C=\{x\in M:f(x)=x\}$ of the
system is the set of fixed points of $f.$ If $g$ is the identity map, the
equilibrium problem can be treated via the so-called Coincidence
Problem\textit{ }\cite[VI.14]{Bredon}, \cite[Chapter 7]{Vick},
\cite{Gorn-book}: \textquotedblleft Given two maps $f,g:N\rightarrow M$
between two $n$-dimensional manifolds, what can be said about the coincidence
set $C$ of all $x$ such that $f(x)=g(x)$?\textquotedblright\ The famous
Lefschetz coincidence theorem states that if the Lefschetz number
$\lambda_{fg}$ is not equal to zero then there is at least one coincidence,
i.e., $C\neq\varnothing.$ Using this and other invariants one can find out
whether a dynamical system has an equilibrium or a periodic point.

In case of a \textit{controlled} dynamical system, the next state $f(x,u)$
depends not only on the current one, $x\in M,$ but also on the \textit{input},
$u\in U.$ Indeed, a discrete time control system is given by the space of
inputs $U$, the space of states $M$, the "state-input" space $N=M\times U,$ a
map $f:N=M\times U\rightarrow M,$ and the projection $g:N=M\times U\rightarrow
M$ (in general $N$ is a fiber bundle and $g$ is the bundle projection). Just
as above, the equilibrium set of the system $C=\{x\in M:f(x,u)=x\}$ is the
coincidence set of the pair $(f,g).$ However since the dimensions of $N$ and
$M$ are not equal anymore, the Lefschetz \textit{number} cannot be defined in
the same fashion. It has to be replaced with the Lefschetz
\textit{homomorphism} \cite{Sav1} which does a better job at detecting
coincidences. Besides equilibria controllability is another application of the
coincidence theory approach. Controllability is treated by proving
surjectivity of certain maps as a map is surjective if it has a coincidence
with any constant map.

The state space $M$ is often a closed manifold when it appears as a
configuration space in robotics. For example, $M=\mathbf{T}^{n}=(\mathbf{S}%
^{1})^{n},$ where $\mathbf{T}^{n}$ is the $n$ dimensional torus, is the
configuration space, i.e., the set of all possible positions, of a robotic arm
with $n$ revolving joints \cite[p. 1]{NvdS}; or $M=\mathbf{R}^{3}\times SO(3)$
is the configuration space of a rigid body in space \cite[Chapter 2]{Latombe}.
Typically, we have $N=M\times U.$ However nontrivial bundles are also common.
For example, consider a spherical pendulum with a gas jet control which is
always directed in the tangent space. Then its state space is $M=\mathbf{S}%
^{2}$, the $2$-sphere, while the state-input space $N$ is the tangent bundle
$T\mathbf{S}^{2}$ of $\mathbf{S}^{2},$ which is an $\mathbf{R}^{2}$-bundle not
isomorphic to $M\times\mathbf{R}^{2}$ \cite[p. 17]{NvdS}. In spite of the
abundance of such examples \cite{BL}, \cite[Chapter 2]{Latombe}, \cite[p.
1]{NvdS}, topological techniques have not thus far found broad applications in
control theory. The only two recent examples known to the author are
Jonckheere \cite{Jon} and Kappos \cite{Kappos-Conley} - \cite{Kappos4}.

A model of a \textquotedblleft plant\textquotedblright\ is a control system as
a triple $(M,N,f)$ of the spaces $M,N$ and the map $f$ as described above.
Since our knowledge of the model is inevitably imprecise, we have to deal with
perturbations of the system. Perturbations may be understood as variations of
the unknown parameters of the system. Therefore, if the system depends
continuously on these parameters, the change of $M,N,$ and $f$ is also
continuous. This means that we have to consider spaces homeomorphic to $M,N$
and maps homotopic to $f.$ An appropriate instrument is homology. Indeed, the
homology of $M,N,f$ remains constant under homeomorphisms and homotopies and
can be rigorously and effectively computed, see Mischaikow \cite{Misch}.

The use of homology provides another benefit. Normally the perturbations of
$f$ are assumed to be \textquotedblleft small\textquotedblright\ (in
particular, this is the basis of the notion of structural stability). However
unless actual estimates are available, we don't know how \textquotedblleft
small\textquotedblright\ these perturbations are in real life. Therefore in
order to take into account the \textquotedblleft worst possible
scenario\textquotedblright\ we should consider \textit{arbitrary homotopies}
of $f.$ In this paper Lefschetz theory is applied to study existence of
equilibria and controllability for systems determined by maps homotopic to $f$
(Theorems \ref{SuffEqDis}, \ref{SuffControl}).

The secondary objective of this paper is to study robustness of some of these
properties under \textquotedblleft small\textquotedblright\ perturbations
because sometimes they may produce dramatic changes in the properties of the
system. We consider situations when this change is the loss of an equilibrium
(Theorem \ref{NecEqDis}) or the loss of controllability (Theorem
\ref{NecContr}) of the system.

The paper is organized as follows. Some preliminaries from algebraic topology
are outlined in the Section \ref{Prel}. In Section
\ref{classical Lefschetz theory} we review the classical theory of Lefschetz
numbers and show its inadequacy for control theory. In Section
\ref{Lefschetz homomorphism} we consider the necessary generalization, the
Lefschetz homomorphism, of the Lefschetz number and state several relevant
results about existence of coincidences. In Section
\ref{Removability of coincidences} we state some results about removability of
coincidences. In Section \ref{Equilibria of discrete systems} we provide
sufficient conditions of existence of equilibria of a discrete system and
their robustness. In Section \ref{Controllability} we provide sufficient
conditions of controllability of a discrete system. In Section
\ref{Controllability of continuous systems} we discuss how our coincidence
results can be applied to existence of equilibria and controllability of
continuous time control systems. Notions of control theory are defined as
needed, for details see \cite{NvdS}, \cite{Sastry}, \cite{Sontag}.

The author thanks the referees for their comments that have helped
significantly improve this paper.

\section{Preliminaries from algebraic topology.\label{Prel}}

In this paper we mostly follow Bredon \cite{Bredon}. Suppose $N$ is a
topological space and $A\subset N$ is a subspace. The singular
\textit{homology group} $H_{k}(N,A)$ of $N$ relative to $A$ over $\mathbf{Q}$
or any other field is defined as follows. If $\Delta_{k}$ is the standard
$k$-simplex, $k=0,1,2...$, any map $\sigma:\Delta_{k}\rightarrow N$ is called
a singular $k$-simplex. We let $C_{k}(N,A)$ be the vector space over
$\mathbf{Q}$ generated by all singular $k$-simplices of $N$ whose images are
not completely in $A.$ Then the boundary operator $\partial_{k}:C_{k}%
(N,A)\rightarrow C_{k-1}(N,A)$ is defined in the natural way. Next we let
$H_{k}(N,A)=\ker\partial_{k}/\operatorname{Im}\partial_{k}.$ Further, let
$C^{k}(N,A)$ be the dual of $C_{k}(N,A),$ i.e., the vector space of all linear
maps from $C_{k}(N,A)$ to $\mathbf{Q.}$ Then $\partial_{k}$ generates the
coboundary operator $\partial^{k}:C^{k}(N,A)\rightarrow C^{k+1}(N,A).$ Next we
let $H^{k}(N,A)=\ker\partial^{k}/\operatorname{Im}\partial^{k}$ be the
\textit{cohomology group} of $N$ relative to $A.$ Also $H_{k}(N)=H_{k}%
(N,\varnothing),$ $H^{k}(N)=H^{k}(N,\varnothing).$ The homology and cohomology
groups $H_{k}(N,A;G),$ $H^{k}(N,A;G)$ over any group $G$ can be defined in a
similar fashion.

Homology and cohomology groups over fields are vector spaces with the
following properties. The Betti numbers, $\dim H_{k}(N),$for $k=0,1,2,$ are
the numbers of path components, "tunnels", and "voids" of $N,$ respectively.
In case of a path connected $N,$ the generators of $H_{0}(N)=H^{0}%
(N)=\mathbf{Q}$ are denoted by $1.$ If $N$ is contractible, it is
\textit{acyclic}, i.e., $H_{k}(N)=H^{k}(N)=0$ for $k>0.$ If $N$ is an
$n$-dimensional simplicial complex, $H_{k}(N)=0$ for all $k>n.$ If $M$ is a
compact connected orientable $n$-dimensional manifold with boundary $\partial
M$ then $H_{n}(M,\partial M)=H^{n}(M,\partial M)=\mathbf{Q.}$ These two groups
are assumed to be generated by the \textit{fundamental classes} $O_{M}$ and
$\overline{O}_{M}$ of $M,$ respectively. Moreover, there is $D_{M}%
:H^{k}(M,\partial M)\rightarrow H_{n-k}(M),$ the \textit{Poincar\'{e} duality}
isomorphism given by the cap product with the fundamental class $O_{M}$. The
\textit{cap product} is the homomorphism $\frown:H^{p}(N,A)\otimes
H_{m}(N,A)\rightarrow H_{m-k}(N)$ given by $x\frown a=(1\times x)\Delta a$,
where $\Delta$ is a diagonal approximation. Then $a\in H_{k}(N,A)$ and $x\in
H^{k}(N,A)$ are called \textit{dual} if $x\frown a=\left\langle
x,a\right\rangle =x(a)=1.$ In particular, $O_{M}$ and $\overline{O}_{M}$ are
dual. By the K\"{u}nneth Theorem, $H_{k}(M\times U)=\sum_{i+j=k}%
H_{i}(M)\otimes H_{j}(U),$ $k=0,1,2,...$.

Suppose $B$ is a subspace of the topological space $B$ and $f:N\rightarrow M$
is a map, then $f:(N,A)\rightarrow(M,B)$ is a \textit{map of pairs} if
$f(A)\subset B.$ In this case the natural homomorphism from $C_{k}(N,A)$ to
$C_{k}(M,B)$ is generated by $f$. This homomorphism generates $f_{\ast}%
:H_{k}(N,A)\rightarrow H_{k}(M,B),$ the \textit{homology homomorphism }of
$f$\textit{, }and\textit{ }$f^{\ast}:H^{k}(M,B)\rightarrow H^{k}(N,A),$ the
\textit{cohomology homomorphism }of $f$. Two maps $f,g:(N,A)\rightarrow(M,B)$
are called \textit{homotopic,} $f\sim g,$ if $f$ can be continuously
\textquotedblleft deformed\textquotedblright\ into $g$, i.e., there is a map
$F:[0,1]\times(N,A)\rightarrow(M,B)$ such that $F(0,\cdot)=f$ and
$F(1,\cdot)=g.$ An $\varepsilon$-\textit{homotopy} is one satisfying
$d(F(t,x),F(0,x))<\varepsilon$ for all $x,t.$ If $f$ and $g$ are homotopic
then $f_{\ast}=g_{\ast}.$ If $f$ is homotopic to a constant map then $f_{\ast
}$ is trivial, i.e., $f_{\ast}:H_{k}(N)\rightarrow H_{k}(M)$ is zero for
$k=1,2,...$. The $k$th homotopy group $\pi_{k}(N)$ of $N$ is the group of
homotopy classes of maps of $k$-spheres to $N$.

\section{Review of Lefschetz theory. \label{classical Lefschetz theory}}

In this section $M$ and $N$ are orientable compact connected manifolds with
boundaries $\partial M,\partial N,$ and $\dim M=\dim N=n$.

Consider the\textbf{\ }Fixed Point Problem: \textquotedblleft If
$f:M\rightarrow M$ is a map, what can be said about the set of points $x\in M$
such that $f(x)=x$?\textquotedblright\ Applications of fixed point theorems
(Kakutani, Banach, etc.) to control problems are abundant, \cite{Balach},
\cite{Carmichael}, \cite{CNZ}, \cite{Id}, \cite{Klamka}, \cite{Nistri}.
However the methods we suggest in this paper go far beyond those.

One may associate to $f$ an integer $\lambda_{f}$ called the Lefschetz number
\cite{Brown}. It detects fixed points and is computable by a simple formula:%
\[
\lambda_{f}=\sum_{n}(-1)^{n}Trace(f_{\ast n}),
\]
where $f_{\ast n}:H_{n}(M)\rightarrow H_{n}(M)$ is induced by $f.$ The
\textit{Lefschetz fixed point theorem} states that if $\lambda_{f}\neq0,$ then
$f$ has a fixed point.

The Coincidence Problem is concerned with a similar question about two maps
$f,g:N\rightarrow M$ and their \textit{coincidences }$x\in N:f(x)=g(x).$ One
of the main tools is the \textit{Lefschetz coincidence number} $\lambda_{fg}$
defined similarly to $\lambda_{f}$ as the alternating sum of traces of a
certain endomorphism on the homology group of $M.$ Algebraically, if
$h:E_{\ast}\rightarrow E_{\ast}$ is a degree $0$ endomorphism of a finitely
dimensional graded vector space $E_{\ast}=\{E_{k}\},h_{k}:E_{k}\rightarrow
E_{k},$ then its \textit{Lefschetz number }is $L(h)=\sum_{k}(-1)^{k}%
Trace(h_{k}).$ To apply this formula we let $E_{\ast}=H_{\ast}(M)$, then the
Lefschetz number is defined as $\lambda_{fg}=L(g_{\ast}D_{N}f^{\ast}D_{M}%
^{-1}),$ where $D_{M}:H^{k}(M,\partial M)\rightarrow H_{n-k}(M),$ $D_{N}%
:H^{k}(N,\partial N)\rightarrow H_{n-k}(N)$ are the Poincar\'{e} duality
isomorphisms. Observe that for $f^{\ast}:H^{k}(M,\partial M)\rightarrow
H^{k}(N,\partial N)$ to be well defined the map $f$ has to be boundary
preserving, $f:(N,\partial N)\rightarrow(M,\partial M).$

A \textit{Lefschetz type coincidence theorem} states that if $\lambda_{fg}%
\neq0$ then $f,g$ (and any pair homotopic to them) have a coincidence$.$ The
converse is false in general. When $\lambda_{fg}=0,$ the maps $f,g$ may have
coincidences but under certain circumstances they can be removed by homotopies
of $f,g$ (see Section \ref{Removability of coincidences}).

Until recently such theorems have been mostly considered in the following two
settings. Case 1: \cite[VI.14]{Bredon}, \cite[Chapter 7]{Vick} $f:N\rightarrow
M$ is a map between two $n$-manifolds as above. This way one can apply the
Lefschetz theorem to detecting equilibria of dynamical systems but not of an
even simplest control system because the dimensions of $N$ and $M$ have to be
equal. Case 2: \cite{Gorn-book} $f:N\rightarrow M$\ is a map from an arbitrary
topological space to an open subset of $\mathbf{R}^{n}$\ and all fibers
$f^{-1}(y)$\ are acyclic, i.e., $H_{k}(f^{-1}(y))=0$\ for\ $k=1,2,...$. Here
the dimensions are also equal in the sense that $H_{\ast}(N)=H_{\ast}(M)$
(Vietoris Theorem). Thus neither case is broad enough to cover all possible
control systems.

As an example from dynamics, one can consider the problem of existence of
closed orbits of a flow given by a map $f:[0,\infty)\times M\rightarrow M,$
i.e., the initial position is $f(0,x)=x$ and $f(t,x)$ gives the position at
time $t$. Closed orbits correspond to coincidences of $f$ and the projection
$p:[0,\infty)\times M\rightarrow M$\textit{. }More generally one considers
$f:X\times M\rightarrow M,$ where $X$ is a topological space. This situation
was studied in \cite{Knill}, \cite{GN2}, \cite{GNO}, \cite{Geoghegan} under
the name \textquotedblleft parametrized fixed point theory\textquotedblright.
These results can be applied to detection of equilibria (Section
\ref{Equilibria of discrete systems}), but the setting is not general enough
to study controllability (Section \ref{Controllability}). The author
\cite{Sav}, \cite{Sav1} extended some of the results of \cite{GNO} to the
general case of two arbitrary maps $f,g:N\rightarrow M$ from an arbitrary
topological space to an orientable compact manifold. The content of these
papers is briefly outlined in the next section.

\section{Detecting coincidences. \label{Lefschetz homomorphism}}

In this section $N$ is an arbitrary topological space, $A\subset N$, $M$ is an
orientable compact connected manifold with boundary $\partial M$, $\dim M=n$,
and $f:(N,A)\rightarrow(M,\partial M),$ $g:N\rightarrow M$ are maps.

The generalization of the Lefschetz number is based on the fact that since the
graded vector space $E=H_{\ast}(M)$ is equipped with the cap product
$\frown:E^{\ast}\otimes E_{\ast}\rightarrow E_{\ast}$, one can define the
Lefschetz class\ $L(h)\in E_{\ast}$ of a graded endomorphism $h$ given by
$h_{k}:E_{k}\rightarrow E_{k+m}$ of any degree $m$ not just of degree $0$ as
in the classical case.

\begin{definition}
\label{Lef-formula}\cite[Proposition 2.2]{Sav1} If $h:H_{k}(M)\rightarrow
H_{k+m}(M),$ $k=0,1,2...,$ is a graded homomorphism of degree $m$ then the
\textit{Lefschetz class}\ is defined as
\[
L(h)=\sum_{k}(-1)^{k(k+m)}\sum_{j}x_{j}^{k}\frown h(a_{j}^{k}),
\]
where $\{a_{1}^{k},...,a_{m_{k}}^{k}\}$ is a basis for $H_{k}(M)$ and
$\{x_{1}^{k},...,x_{m_{k}}^{k}\}$ the corresponding dual basis for $H^{k}(M).$
\end{definition}

For a given $z\in H_{s}(N,A),$ suppose $h_{fg}^{z}$ is defined as the
composition%
\[
H_{i}(M)^{\underrightarrow{\quad D_{M}^{-1}\quad}}H^{n-i}(M,\partial
M)^{\underrightarrow{\quad f^{\ast}\quad}}H^{n-i}(N,A)^{\underrightarrow
{\quad\frown z\quad}}H_{s-n+i}(N)^{\underrightarrow{\quad g_{\ast}\quad}%
}H_{s-n+i}(M),
\]
i.e.,%
\[
h_{fg}^{z}(x)=g_{\ast}((f^{\ast}D_{M}^{-1}(x))\frown z)).
\]
Its degree is $m=s-n.$

\begin{definition}
The \textit{Lefschetz homomorphism} $\Lambda_{fg}:H_{k}(N,A)\rightarrow
H_{k-n}(M),$ $k=0,1,...,$ is defined by
\[
\Lambda_{fg}(z)=L(h_{fg}^{z}).
\]

\end{definition}

\begin{proposition}
\label{Traces}If the degree $m$ of $h$ is zero, $L(h)=\sum_{k}(-1)^{k}%
Trace(h_{k})$.
\end{proposition}

In particular, the degree of the homomorphism $h_{fg}^{z}$ is zero if $z\in
H_{n}(N,A).$ If, moreover, $N$ is a orientable compact connected manifold of
dimension $n,$ we have $H_{n}(N,\partial N)=\mathbf{Q.}$ It is generated by
the fundamental class $O_{N}\in H_{n}(N,\partial N)$ of $N.$ Since
$D_{N}(x)=x\frown O_{N},$ we recover the classical \textit{Lefschetz number},
$\lambda_{fg}=\Lambda_{fg}(O_{N}).$

\begin{theorem}
\label{Lef-type}\cite[Theorem 6.1]{Sav1} \textbf{(Existence of coincidences)}
If\textit{ }$\Lambda_{fg}\neq0$\textit{\ then any pair of maps }$f^{\prime
},g^{\prime}$\textit{\ homotopic to }$f,g$ has a coincidence.
\end{theorem}

Especially important for control theory are the following corollaries. They
are applied to existence of equilibria (Section
\ref{Equilibria of discrete systems}), controllability (Sections
\ref{Controllability} and \ref{Controllability of continuous systems}).
Observe that the second corollary is about a map of pairs and the first is not.

\begin{corollary}
\label{Proj}\textbf{(Existence of fixed points)} (cf. \cite{GNO}) Let
$g:M\times U\rightarrow M$ be a map. Given $v\in H_{s}(U),$ suppose the
homomorphism $g_{v}:H_{i}(M)\rightarrow H_{i+s}(M)$ of degree $s$ is defined
by
\[
g_{v}(x)=(-1)^{(n-i)s}g_{\ast}(x\otimes v),
\]
$x\in H_{i}(M).$ Then, if
\[
L(g_{v})\neq0\text{ for some}\ v\in H_{s}(U)
\]
then any map $g^{\prime}:M\times U\rightarrow M$ homotopic to $g$ has a fixed
point $x$, $g^{\prime}(x,u)=x$ for some $u.$
\end{corollary}

\begin{proof}
Suppose $(N,A)=(M,\partial M)\times U$ and apply the above theorem to the pair
$p,g,$ where $p$ is the projection $p:(M,\partial M)\times U\rightarrow
(M,\partial M).$ Also according to Corollary 5.7 in \cite{Sav1}, $\Lambda
_{pg}(O_{M}\otimes v)=L(g_{v})$.
\end{proof}

\begin{corollary}
\label{Surj}\textbf{(Sufficient condition of surjectivity)} If
\[
f_{\ast}:H_{n}(N,A)\rightarrow H_{n}(M,\partial M)=\mathbf{Q}\mathit{\ }%
\text{\textit{is nonzero}}%
\]
\textit{ then any }map $f^{\prime}:(N,A)\rightarrow(M,\partial M)$
\textit{homotopic to }$f$ is onto.
\end{corollary}

\begin{proof}
Apply the theorem to the pair $f,c,$ where $c$ is any constant map (as in
Section 5 in \cite{Sav} and Proposition 6.8 in \cite{Sav1}).
\end{proof}

In case of manifolds of equal dimensions the condition of this corollary is
equivalent to nonvanishing of the topological degree \cite[p. 186]{Bredon} of
$f$.

\section{Removing coincidences.\label{Removability of coincidences}}

In this section $M$ is a compact orientable connected manifold, $\dim M=n,$
$N$ is a manifold, $A$ is a closed subset of $N$, $f:(N,A)\rightarrow
(M,\partial M),$ $g:N\rightarrow M$ are maps.

When $\dim N=\dim M=n>2,$ the vanishing of the Lefschetz number $\lambda_{fg}$
implies that the coincidence set can be removed by homotopies of $f,g$
\cite{BS}. If $\dim N=n+m,m>0,$ this is no longer true even if $\lambda_{fg}$
is replaced with $\Lambda_{fg}$. Some progress has been made for $m=$ $1.$ In
this case the secondary obstruction to the removability of a coincidence set
was considered in \cite{Fuller}, \cite{DG}, \cite{Jez}. These results can be
used to study removability of equilibria when the dimension of the input space
is $1$. However the conditions on $f$ and $g$ are hard to verify. Necessary
conditions of the global removability for arbitrary $m$ were considered in
\cite[Section 5]{GJW} with $N$ a torus and $M$ a nilmanifold. For some $m>1,$
a partial converse of the Theorem \ref{Lef-type} is provided by the author
\cite{Sav2}. A version of this theorem is given below.

Let $F$ be a closed subset of $N.$ We say that $F$ satisfies condition (*)
(condition (A) in \cite[Section 3]{Sav2}) if there are arbitrarily small
neighborhoods $W$ and $V$ of $F$ such that $V\subset\overline{V}\subset
W\subset N\backslash A$ and
\[
\text{(*) \ }H^{k+1}(W,W\backslash V;\pi_{k}(\mathbf{S}^{n-1}))=0\text{ for
}k\geq n+1.
\]
In particular this condition is satisfied if $F$ is contractible.

\begin{theorem}
\label{Remov}\textbf{(Removability of coincidences) }Suppose condition (*) is
satisfied for $F,$ an isolated subset of the coincidence set of $f,g$. Suppose
also $f(F)=g(F)=\{x\},$ $x\in M\backslash\partial M.$ Then, if
\[
\Lambda_{fg}(z)=\mathit{\ }\sum_{k}(-1)^{k}Trace(h_{fg}^{z})\in H_{0}%
(M)=\mathbf{Q}\text{ \textit{is zero for all}}\mathit{\ }z\in H_{n}(N,A)
\]
\textit{then there is a homotopy of }$f$ (or $g)$ to a map $f^{\prime}$ (or
$g^{\prime})$\textit{ such that the new pair has no coincidences}. The
homotopy can be chosen \textit{arbitrarily small and constant on the
compliment of a neighborhood of }$F.$
\end{theorem}

\begin{proof}
Since $f(F)=\{x\},$ $x\in M\backslash\partial M,$ and $f(A)\subset\partial M$,
we have $F\cap A=\varnothing.$ Therefore the map $(f,g):(N,A)\rightarrow
M^{\times}=(M\times M,M\times M\backslash d(M))$, where $d(M)$ is the diagonal
of $M\times M,$ is well defined. Let $I_{fg}^{N}(\tau)\in H^{n}(N,A)$ be the
cohomology coincidence index \cite[Section 2]{Sav2}. Here $I_{fg}%
^{N}=(f,g)^{\ast}:H^{n}(M^{\times})\rightarrow H^{n}(N,A)$, $M^{\times},$
$\tau$ is the generator of $H^{n}(M^{\times})=\mathbf{Q}$. Let $I_{fg}$ be the
homology coincidence homomorphism defined by $I_{fg}=(f,g)_{\ast}%
:H_{k}(N,A)\rightarrow H_{k}(M^{\times}).$ By Theorem 6.1 in \cite{Sav1},
$\Lambda_{fg}(z)=\pi_{\ast}(\tau\frown I_{fg}(z)),$ where $\pi:M\times
M\rightarrow M$ is the projection on the first factor. Then, for any $z\in
H_{n}(N,A)$%
\begin{align*}
\Lambda_{fg}(z)  &  =\pi_{\ast}(\tau\frown(f,g)_{\ast}(z))=\pi_{\ast
}(f,g)_{\ast}((f,g)^{\ast}(\tau)\frown z)\\
&  =\left\langle (f,g)^{\ast}(\tau),z\right\rangle =\left\langle I_{fg}%
^{N}(\tau),z\right\rangle .
\end{align*}
Therefore $\Lambda_{fg}(z)=0$ for all $z\in H_{n}(N,A)$ if and only if
$I_{fg}^{N}(\tau)=0.$ The trace formula for $\Lambda_{fg}$ comes from
Proposition \ref{Traces}.

Let $G$ be an open neighborhood of $x.$ Choose neighborhoods $W$ and $V$ of
$F=f^{-1}(x)$ as in condition (*) so that $f(W)\subset G$.\ Now we proceed as
in the proof of Theorem 2 in \cite{Sav2} to show that the coincidence subset
$F$ can be removed by a homotopy of $f$ or $g$ relative to $N\backslash W$
provided the local cohomology index $I_{fg}^{W}(\tau)$ vanishes. This index is
defined as the one above: $I_{fg}^{W}=(f,g)^{\ast}:H^{n}(M^{\times
})\rightarrow H^{n}(W,W\backslash V).$ Now if $i:W\rightarrow N$ is the
inclusion, then $I_{fg}^{W}(\tau)=i^{\ast}I_{fg}^{N}(\tau)=0.$
\end{proof}

Condition (*) ensures that only the primary obstruction to removability, i.e.,
the Lefschetz number, may be nonzero. Further investigation of necessary
conditions of removability of coincidences will require consideration of
higher order obstructions. The case when $f(F)$ is not a single point is best
addressed in the context of Nielsen theory via Wecken type theorems
\cite{Sav4}. In general, the homotopy of $f$ is not local.

Especially important for control theory are the following corollaries. They
are applied to disappearance under perturbations of equilibria (Sections
\ref{Equilibria of discrete systems} and
\ref{Controllability of continuous systems}) and controllability (Section
\ref{Controllability}).

\begin{corollary}
\label{NecParam}\textbf{(Removability of fixed points)} Suppose $U$ is a
manifold and suppose $a\in M\backslash\partial M$ is an isolated fixed point
of a map $g:M\times U\rightarrow M$, i.e., $g(a,u)=a$ for some $u\in U.$
Suppose condition (*) is satisfied for $F=\{u\in U:g(a,u)=a\}$. Then, if
\[
L(g_{1})=\sum_{k}(-1)^{k}Trace(\overline{g}_{\ast k})=0,
\]
where $\overline{g}(\cdot)=g(\cdot,u_{0}),$ then \textit{there is a homotopy
of }$g$ to a map $g^{\prime}$\textit{ such that} $g^{\prime}$ has no fixed
points$.$ The homotopy can be chosen \textit{arbitrarily small and constant on
the compliment of a neighborhood of }$\{a\}\times F.$
\end{corollary}

\begin{proof}
Suppose $(N,A)=(M,\partial M)\times U$ and apply the theorem to the pair
$p,g,$ where $p$ is the projection $p:(M,\partial M)\times U\rightarrow
(M,\partial M).$ Next we use the formula $\Lambda_{pg}(O_{M}\otimes
v)=L(g_{v})$ from Corollary 5.7 in \cite{Sav1}. Finally, we observe that the
coincidence set $\{a\}\times F$ satisfies condition (*).
\end{proof}

\begin{corollary}
\label{NecSurj}\textbf{(Removability of images)} Suppose that there is a fiber
$F=f^{-1}(x_{0})$ of $f$ satisfying condition (*).\ Then, if
\[
f_{\ast}:H_{n}(N,A)\rightarrow H_{n}(M,\partial M)=\mathbf{Q}\mathit{\ }%
\text{\textit{is zero}}%
\]
\textit{then there is a homotopy of }$f$ to a map $f^{\prime}$ which is not
onto; specifically, $x_{0}\notin f^{\prime}(N)$. The homotopy can be chosen
\textit{arbitrarily small and constant on the compliment of a neighborhood of
}$F.$
\end{corollary}

\begin{proof}
Suppose $c(x,u)=x_{0}$ is the constant map. Next, $f_{\ast}=0$ if and only if
$\Lambda_{fc}(z)=0$ for all $z\in H_{n}(N,A)$ (see Section 5 in \cite{Sav}).
Now apply the theorem to the pair $f,c$ (cf. Theorem 3 in \cite{Sav2}).
Observe also that $f^{-1}(x_{0})$ is the coincidence set of the pair.
\end{proof}

A closed submanifold $F$ of $N$ satisfies condition (*) if one of the
following\ three conditions holds \cite[Section 4]{Sav2}:

\begin{itemize}
\item \textit{(a1) }$M$\textit{ is a surface, i.e., }$n=2;$\textit{ or}

\item \textit{(a2) }$F$\textit{ is acyclic, i.e., }$H_{k}(F)=0$%
\textit{\ for\ }$k=1,2,...$\textit{; \ or}

\item \textit{(a3) every component of }$F$\textit{ is a homology }%
$m$\textit{-sphere, i.e., }$H_{k}(F)=0$\textit{\ for\ }$k\neq0,m,$\textit{ for
the following values of }$m$\textit{ and }$n$\textit{:}
\end{itemize}

\textit{\qquad(1) }$m=4$\textit{ and }$n\geq6;$

\textit{\qquad(2) }$m=5$\textit{ and }$n\geq7;$

\textit{\qquad(3) }$m=12$\textit{ and }$n=7,8,9,$ or $n\geq14$\textit{.}

\section{Existence of equilibria.\label{Equilibria of discrete systems}}

In this section $M$\textit{ }is a compact orientable connected manifold with
boundary $\partial M$, $\dim M=n$, $U$ is a topological space, $U^{\prime
}\subset U$.

A \textit{discrete time control system} $D_{g}$ is given by a map $g:M\times
U\rightarrow M$, with $U$ the \textit{space of inputs}, $M$ the \textit{space
of states} of the system.

We say that $D_{g^{\prime}}$ is a \textit{perturbation} of $D_{g}$ if
$g^{\prime}$ homotopic to $g.$ To justify this definition recall that the
system $D_{g^{\prime}}$ is normally called a perturbation of $D_{g}$ if
$g^{\prime}$ satisfies $\rho(g(z),g^{\prime}(z))<\varepsilon$ for all $z\in
M\times U$, where $\rho$ is the distance on $M,$ for some $\varepsilon>0$.
However, if $M\ $is a manifold, the above condition for a small enough
$\varepsilon$ implies that $g$ and $g^{\prime}$ are homotopic.

As before suppose $\{a_{1}^{k},...,a_{m_{k}}^{k}\}$ is a basis for $H_{k}(M)$
and $\{x_{1}^{k},...,x_{m_{k}}^{k}\}$ the corresponding dual basis for
$H^{k}(M).$

\begin{theorem}
\textbf{\label{SuffEqDis}(Existence of equilibria) }If
\[
L(g_{v})=(-1)^{ns}\sum_{k}(-1)^{k}\sum_{j}x_{j}^{k}\frown g_{\ast}(a_{j}%
^{k}\otimes v)\neq0\text{ for some }v\in H_{s}(U)
\]
then every perturbation of the discrete time system $D_{g}$ has an equilibrium.
\end{theorem}

\begin{proof}
In light of Corollary \ref{Proj} we only need to show that the above formula
for the Lefschetz number $L(g_{v})$ of $g_{v}(x)=(-1)^{(n-i)s}g_{\ast
}(x\otimes v),$ $x\in H_{i}(M),$ is true. The degree of $g_{v}$ is $s,$
$a_{j}^{k}\in H_{k}(M).$ Therefore we substitute $m=s$ and $i=k$ in Definition
\ref{Lef-formula}:%
\begin{align*}
L(g_{v})  &  =\sum_{k}(-1)^{k(k+s)}\sum_{j}x_{j}^{k}\frown(-1)^{(n-k)s}%
g_{\ast}(a_{j}^{k}\otimes v)\\
&  =\sum_{k}(-1)^{k^{2}+ns}\sum_{j}x_{j}^{k}\frown g_{\ast}(a_{j}^{k}\otimes
v),
\end{align*}
and the formula follows.
\end{proof}

\begin{corollary}
\textbf{\label{SuffEqDisCor}}Suppose $M=\mathbf{S}^{n},$ and suppose one of
the following conditions is satisfied:

(1) $g_{\ast}(d\otimes1)\neq(-1)^{n+1}d,$ where $d$ is a generator of
$H_{n}(\mathbf{S}^{n});$ or

(2) $g_{\ast}(1\otimes v)\neq0$ for some $v\in H_{n}(U).$

Then every perturbation of the discrete time system $D_{g}$ has an equilibrium.
\end{corollary}

\begin{proof}
Recall that $a_{j}^{k}\otimes v\in H_{k+s}(M).$ Since $H_{k}(\mathbf{S}%
^{n})=0$ for all $k\neq0,n,$ and $g(x,u)\in\mathbf{S}^{n}$, we have $g_{\ast
}(a_{j}^{k}\otimes v)=0$ except for the following two cases. (1) $v=1\in
H_{0}(U),s=0,$ then $k=0,a_{j}^{0}=x_{j}^{0}=1,$ or $k=n,a_{j}^{n}=d,x_{j}%
^{n}=\overline{d},$ or (2) $v\in H_{n}(U),s=n,$ then $k=0,a_{j}^{0}=x_{j}%
^{0}=1.$ Here $\overline{d}$ is the dual of $d,$ $\overline{d}\frown d=1.$
Thus we have%
\begin{align*}
\text{(1) }L(g_{1})  &  =(-1)^{n0}(1\frown g_{\ast}(1\otimes1)+(-1)^{n}%
\overline{d}\frown g_{\ast}(d\otimes1))\\
&  =1+(-1)^{n}\overline{d}\frown g_{\ast}(d\otimes1);\\
\text{(2) }L(g_{v})  &  =(-1)^{nn}(1\frown g_{\ast}(1\otimes v))\\
&  =(-1)^{n}g_{\ast}(1\otimes v).
\end{align*}
Now, if either $L(g_{1})$ or $L(g_{v})$ is nonzero, then $D_{g}$ has an
equilibrium by the above theorem.
\end{proof}

Condition (1) means that the degree of $\overline{g}(\cdot)=g(\cdot
,u_{0}):\mathbf{S}^{n}\rightarrow\mathbf{S}^{n}$ is not equal to $(-1)^{n+1}.$

If $g:M\times M\rightarrow M$ is the multiplication of a compact Lie group,
then $D_{f}$ has a equilibrium \cite[Example 2.3]{GNO}. For more examples, see
\cite{GNO}, \cite{Sav}, \cite{Sav1}.

In the control setting Corollary \ref{NecParam} reads as follows.

\begin{theorem}
\textbf{\label{NecEqDis}(Robustness of equilibria) }Suppose $U$ is a manifold
and suppose $a\in M\backslash\partial M$ is an isolated equilibrium of
$D_{g}.$ Suppose condition (*) is satisfied for $F=\{u\in U:g(a,u)=a\}$. Then,
if
\[
L(g_{1})=\sum_{k}(-1)^{k}Trace(\overline{g}_{\ast k})=0,
\]
where $\overline{g}(\cdot)=g(\cdot,u_{0}),$ then this equilibrium can be
removed by an \textit{arbitrarily} small perturbation restricted to a
neighborhood of $F$.
\end{theorem}

In particular if $M=\mathbf{S}^{n},$ then $L(g_{1})=1+(-1)^{n}\deg\overline
{g}.$ Therefore in light of condition (1) of the above corollary we conclude
that $D_{g}$ has a robust equilibrium if and only if the degree of
$\overline{g}$ is not equal to $(-1)^{n+1}$.

\section{Controllability.\label{Controllability}}

In this section $M$ is a compact orientable connected manifold with boundary
$\partial M,$ $\dim M=n$, $U$ is a topological space.

Suppose a discrete system $D_{f}$ is given by $f:M\times U\rightarrow M.$ The
system $D_{f}$ is called \textit{controllable} \cite{Sontag} if any state can
be reached from any other state, i.e., for each pair of states $x,y\in M$
there are inputs $u_{0},...,u_{r}\in U$ such that $x_{1}=f(u_{0}%
,x),x_{2}=f(u_{1},x_{1}),...,y=x_{r}=f(u_{r},x_{r}),$ notation
$x\rightsquigarrow_{f}y.$

Below this notion is generalized in three, nontypical but topologically
appropriate, ways. First, we consider the possibility of an arbitrary state
reached not from any given state but from a state in a particular subset $L$
of $M.$ Second, to define robustness we allow for arbitrary, not necessarily
small, perturbations of $f.$ Third, instead of looking into controllability of
a new, perturbed, system $D_{g},$ where $g$ is a perturbation of $f$, we allow
for consecutive applications of different maps homotopic to $f$.

To apply Corollary \ref{Surj} we need $f$ to be a map of pairs. For this
purpose in this section we make the following assumption about $D_{f}$. If the
initial state lies at the boundary $\partial M$ of $M$ then the next state,
regardless of the input, lies within a certain neighborhood of $\partial M$.
For simplicity we make a simpler, but topologically equivalent, assumption,
$f(\partial M\times U)\subset\partial M$.

Next, let $U^{\prime}$ be the set of controls that take any given state to the
boundary of $M,$ i.e.,
\[
U^{\prime}=\{u\in U:f(x,u)\in\partial M\text{ for all }x\in M\}.
\]
Therefore $f(M\times U^{\prime})\subset\partial M.$ Combining this with the
above assumption we can treat $f$ as a map of pairs, $f:(M,\partial
M)\times(U,U^{\prime})\rightarrow(M,\partial M).$

\begin{definition}
Given $L\subset M,$ let $f^{\prime}:L\times U\rightarrow M$ be the restriction
of $f.$ Then the system is called \textit{strongly robustly controllable from
}$L$ if there is $\varepsilon>0$ such that for any map $f_{0}$ homotopic to
$f^{\prime},$ any maps $f_{1},...,f_{r}$ homotopic to $f,$ and for each $y\in
M$ there is $x\in L$ and there are inputs $u_{0},...,u_{r}\in U$ such that
\[
x_{1}=f_{0}(x,u_{0}),x_{2}=f_{1}(x_{1},u_{1}),...,y=x_{r+1}=f_{r}(x_{r}%
,u_{r}).
\]

\end{definition}

The system is controllable if it controllable from a point.

Observe that $f^{\prime}:(L,L^{\prime})\times(U,U^{\prime})\rightarrow
(M,\partial M)$ is a map of pairs, where $L^{\prime}=L\cap\partial M.$

The following theorem translates the above "reachability" condition into the
language of homology: any element of $H_{n}(M,\partial M)=\mathbf{Q},$ and in
particular the fundamental class $O_{M}$ of $M,$ can be reached from some
$a_{0}\in H_{\ast}(L,L^{\prime})$ by means of $f_{\ast}$.

\begin{theorem}
\label{SuffControl}\textbf{(Sufficient condition of robust controllability)
}Suppose that there are $a_{0}\in H_{p}(L,L^{\prime}),$ $v_{0}\in H_{s_{0}%
}(U,U^{\prime}),...,v_{r}\in H_{s_{r}}(U,U^{\prime})$ such that
\[
a_{1}=f_{\ast}^{\prime}(a_{0}\otimes v_{0}),a_{2}=f_{\ast}(a_{1}\otimes
v_{1}),...,a_{r+1}=f_{\ast}(a_{r}\otimes v_{r})\in H_{n}(M,\partial
M)\backslash\{0\}.
\]
Then the discrete time system $D_{f}$ is strongly \textit{robustly
}controllable from $L$.
\end{theorem}

Here, if $a_{i}\in H_{n_{i}}(M,\partial M),i=0,1,2,...,r,$ then $n_{0}%
=p,n_{1}=p+s_{0},n_{2}=n_{1}+s_{1},...,n_{r}=n_{r-1}+s_{r}=n.$ Thus we have a
sequence of homology classes $a_{0},..,a_{r}$ "climbing" dimensions from $p$
to $n.$

\begin{proof}
The result of consecutive applications of $f$ is defined as a map
$F:(L,L^{\prime})\times(U,U^{\prime})^{r+1}\rightarrow(M,\partial M)$ given
by
\[
F(x,u_{0},...,u_{r})=f(...f(f^{\prime}(x,u_{0}),u_{1}),...,u_{r}),
\]
i.e., it is given by the composition
\begin{align*}
F  &  :(L,L^{\prime})\times(U,U^{\prime})\times...\times(U,U^{\prime
})^{\underrightarrow{~\ \ \ f^{\prime}\times Id\ \ \ ~~}}\\
&  (M,\partial M)\times(U,U^{\prime})\times...\times(U,U^{\prime
})^{\underrightarrow{~\ \ \ f\times Id\ \ \ ~~}}...\text{.}%
\end{align*}
Then $x\rightsquigarrow_{f}F(x,u_{0},...,u_{r}).$ Suppose a map $f_{0}$ is
homotopic to $f^{\prime}$ and maps $f_{1},...,f_{r}$ are homotopic to $f.$ The
result of consecutive applications of $f_{0},...,f_{r}$ is defined as a map
$G:(L,L^{\prime})\times(U,U^{\prime})^{r+1}\rightarrow(M,\partial M)$ given
by
\[
G(x,u_{0},...,u_{r})=f_{r}(...f_{1}(f_{0}(x,u_{0}),u_{1}),...,u_{r}).
\]
Therefore strong robust controllability from $L$ means that $G_{r}:L\times
U^{r+1}\rightarrow M$ is onto. By Corollary \ref{Surj} if
\[
F_{\ast}:H_{n}((L,L^{\prime})\times(U,U^{\prime})\times...\times(U,U^{\prime
}))\rightarrow H_{n}(M,\partial M)=\mathbf{Q}%
\]
is nonzero then every map homotopic to $F$ is onto. Since $G$ is clearly
homotopic to $F,$ all we need to prove is that $F_{\ast}$ is nonzero. By the
K\"{u}nneth theorem $F_{\ast}$ is given by the composition
\begin{align*}
F_{\ast}  &  :H_{\ast}(L,L^{\prime})\otimes H_{\ast}(U,U^{\prime}%
)\otimes...\otimes H_{\ast}(U,U^{\prime})^{\underrightarrow{~\ \ \ f_{\ast
}^{\prime}\otimes Id\ \ \ ~~}}\\
&  H_{\ast}(M,\partial M)\otimes H_{\ast}(U,U^{\prime})\otimes...\otimes
H_{\ast}(U,U^{\prime})^{\underrightarrow{~\ \ \ f_{\ast}\otimes Id\ \ \ ~~}%
}...\text{.}%
\end{align*}
Now the condition of the theorem implies that $f_{\ast}(...f_{\ast}(f_{\ast
}^{\prime}(a_{0}\otimes v_{0})\otimes v_{2})\otimes...\otimes v_{r})\neq0$ for
some $a_{0}\in H_{p}(L,L^{\prime})$ and some $v_{0}\in H_{s_{1}}(U,U^{\prime
}),...,v_{r}\in H_{s_{r}}(U,U^{\prime})$ such that $p+s_{1}+...+s_{r}=n.$
Therefore $F_{\ast}(a_{0}\otimes v_{0}\otimes v_{2}\otimes...\otimes
v_{r})\neq0$.
\end{proof}

Moreover, it is clear that what we have is the \textquotedblleft finite time
reachability\textquotedblright, i.e., every state can be reached in a finite
number, $r+1,$ of steps and that number is common for all states.

The theorem involves multiple iterations of $f_{\ast}$ while it is preferable
to have a condition involving only $f_{\ast}$ itself. Let's consider a case
when this is possible.

Consider first a simple example, $U=\mathbf{S}^{1},$ $U^{\prime}=\varnothing,$
$M=\mathbf{T}^{n}$, and $f:\mathbf{S}^{1}\times\mathbf{T}^{n}\rightarrow
\mathbf{T}^{n}$ is given by $f(u,x_{1},...,x_{n})=(u,x_{1},...,x_{n-1}).$ This
may serve as a model for a robotic arm with $n$ joints where only the first
joint can be controlled directly and the next state of a joint is
\textquotedblleft read\textquotedblright\ from the current state of the
previous joint. The system is obviously controllable. Indeed after $n$
iterations with inputs $u_{1},...,u_{n}$ the system's state is $(u_{n}%
,...,u_{1}).$ Whether the system is robustly controllable is not as obvious.
The affirmative answer is provided by the theorem. Indeed let $L$ be a point,
$p=0.$ Now, if $d$ is the generator of $H_{1}(\mathbf{S}^{1})$, we choose%
\begin{align*}
v_{0}  &  =v_{1}=...=v_{n}=d\in H_{1}(\mathbf{S}^{1}),\text{ and}\\
a_{0}  &  =1\in H_{0}(\mathbf{S}^{1}),\\
a_{1}  &  =d\in H_{1}(\mathbf{S}^{1}),\\
a_{2}  &  =d\otimes d\in H_{2}(\mathbf{S}^{1}\times\mathbf{S}^{1}),\\
&  ...\\
a_{n}  &  =d\otimes...\otimes d\in H_{n}(\mathbf{S}^{1}\times...\times
\mathbf{S}^{1}).
\end{align*}

More generally, suppose the state space $M$ has the product structure,
$M=K_{1}\times...\times K_{s},$ where $K_{i}$ are manifolds of dimensions
$k_{i}.$ Suppose $f=(h_{1},...,h_{s})$, where $h_{i}:U\times M\rightarrow
K_{i}.$ Suppose for $i=1,...,s,$ maps $h_{i}^{a}:K_{i-1}\rightarrow K_{i},$
where $K_{0}=U,$ are given by $h_{i}^{a}(x_{i-1})=h_{i}(a_{0},...,a_{i-2}%
,x_{i-1},a_{i},...,a_{s}).$ If all $h_{i}^{a}$ are onto then the system is
controllable. According to Corollary \ref{Surj} it suffices to require that
all $h_{i\ast}^{a}:H_{k_{i}}(K_{i-1})\rightarrow H_{k_{i}}(K_{i})$ are
nonzero, $i=1,...,s$.

The above theorem can be informally understood as follows. If there are some
submanifolds $M_{1},...,M_{r}$ of $M$ such that $M_{1}=f(L\times U),$
$M_{2}=f(M_{1}\times U),...,M=f(M_{r}\times U)$ then the system is
controllable. Therefore the restrictions $f_{0}:L\times U\rightarrow
M_{1},f_{1}:M_{1}\times U\rightarrow M_{2},...,f_{r}:M_{r}\times U\rightarrow
M,$ of $f$ are onto$.$ This holds provided $f_{i\ast}(O_{M_{i}}\otimes
O_{U})=O_{M_{i+1}},$ where $\dim M_{i}=n_{i},$ $O_{M_{i}}\in H_{n_{i}}(M_{i})$
is the fundamental class of $M_{i},$ $i=0,1,...,r,$ $M_{0}=L.$ Since each
$O_{M_{i}}$ corresponds to $a_{i}=J_{i\ast}(O_{M_{i}})\in H_{n_{i}}(M),$ where
$J_{i}:$ $M_{i}\rightarrow M$ is the inclusion, we arrive at the requirement
of the theorem. The robustness of each of these surjectivity conditions can be
tested by means of Corollary \ref{NecSurj}. As two special cases we have the following.

\begin{theorem}
\label{NecContr}\textbf{(Necessary condition of robust controllability)
}Suppose\textbf{ }$U$ is a manifold and there is a fiber $f^{-1}(x),x\in M,$
of $f$ satisfying condition (*).

(1) Suppose
\[
\text{ }f_{\ast}:H_{n}((M,\partial M)\times(U,U^{\prime}))\rightarrow
H_{n}(M,\partial M)\text{ is zero.}%
\]
Then there is an \textit{arbitrarily} small perturbation of the discrete time
system $D_{f}$ which is not controllable from $M;$ specifically, $x$ is
unreachable from any point.

(2) Suppose $x$ cannot be reached in one step from $L\subset M,$ i.e.,
$x\notin f(L,u).$ Suppose
\[
\text{ }f_{\ast}:H_{n}((M\backslash D,\partial M)\times(U,U^{\prime
}))\rightarrow H_{n}(M,\partial M)\text{ is zero,}%
\]
where $D$ is a small disk around $x.$ Then there is an \textit{arbitrarily}
small perturbation of the discrete time system $D_{f}$ which is not
controllable from $L;$ specifically, $x$ is unreachable from $L$.
\end{theorem}

\section{Continuous systems.\label{Controllability of continuous systems}}

In this section we outline, in less details than above, the possibilities of
applying Lefschetz numbers to continuous systems.

In this section $M$\textit{ }is a compact orientable connected smooth manifold
with boundary $\partial M$, $\dim M=n.$ Let $TM$ be the tangent bundle of $M,$
then $\dim TM=2n$.

A \textit{continuous time control system} $C_{h}$ \cite[p. 16]{NvdS} is
defined as a commutative diagram
\[%
\begin{tabular}
[c]{ccc}%
$Q$ & $^{\underrightarrow{\quad h\quad}}$ & $TM,$\\
$\downarrow^{p}$ & $\swarrow^{\pi_{M}}$ & \\
$M$ &  &
\end{tabular}
\ \ \ \ \ \
\]
where $p:Q\rightarrow M$ is a fiber bundle over $M$ and $\pi_{M}$ is the
projection$.$ In other words we have a parametrized vector field on $M.$ We
say that $C_{k}$ is a \textit{perturbation} of $C_{h}$ if $k$ homotopic to
$h.$

We say that $x\in M$ is an \textit{equilibrium }of this system\textit{ }if
there is $z\in Q$ such that $h(z)=(x,0)\in TM.$ Detecting an equilibrium can
be restated as a coincidence problem. Suppose $i:M\rightarrow TM$ is the
inclusion and $p_{1}:Q\times M\rightarrow Q,$ $p_{2}:Q\times M\rightarrow M$
are the projections. Define the maps $f,g:Q\times M\rightarrow TM$ by
$f=hp_{1},g=ip_{2}.$ Then a coincidence of the pair $f,g$ is an equilibrium of
the system $C_{h}.$ Therefore equilibria can be detected by means of the
coincidence results in Section \ref{Lefschetz homomorphism} and their
robustness can be studied by means of Section
\ref{Removability of coincidences}.

We have a simpler coincidence problem if $M\ $is parallelizable, i.e., $TM$ is
isomorphic to $M\times\mathbf{R}^{n}$. For example, $\mathbf{S}^{1},$
$\mathbf{S}^{3},$ $\mathbf{S}^{7}$ are parallelizable. Let $q:TM\simeq
M\times\mathbf{R}^{n}\rightarrow M$ be the projection. Then a coincidence of
the pair $qh,p$ is an equilibrium of the system $C_{h}$ and we can use Theorem
\ref{Lef-type} to detect equilibria and Theorem \ref{Remov} to study their
robustness. In fact $D_{qh}$ is a discrete system associated with the
continuous system $C_{h}.$ In particular, when $Q=M\times U,$ the results of
Sections \ref{Equilibria of discrete systems} and \ref{Controllability} can be
applied to study equilibria and controllability of $C_{h}.$

When $M\ $is not parallelizable, a discrete system $D_{f}$ associated to the
continuous system $C_{h}$ may be constructed as follows.

Let $\mathcal{A}$ be the topological space of \textit{admissible controls
}associated with $C_{h}$\textit{, }i.e., a set of maps\ $z:[0,d]\rightarrow
Q,$ for all $d\in\mathbf{R}.$ A map $c_{z}:[0,d]\rightarrow M$ is called a
\textit{trajectory} of the control system if there exists a control\ $z\in
\mathcal{A}$ satisfying: $pz=c_{z}$ and $\dfrac{d}{dt}c_{z}=hz.$

Suppose $Q=M\times U,$ where $U$ is the topological space of all possible
inputs, and $p:Q=M\times U\rightarrow M$ is the projection, then $\mathcal{A}$
is the set of pairs $(c,p),$ where $c:[0,d]\rightarrow M$ is a trajectory and
$p:[0,d]\rightarrow U$ is an input$.$ To simplify things even further we
consider only constant inputs. First we assume that the system $C_{h}$
satisfies the following existence and uniqueness property \cite{Sontag}: for
every $x\in M$ and any \textit{constant} input $p(t)=u\in U$ there is a unique
trajectory $c$ such that $c(0)=x$ and $(c,u)\in\mathcal{A}$. Then the
following end point map $f_{d}:M\times U\rightarrow M$ is well defined. Let
$f_{d}(x,u)=c(d),$ where $c:[0,d]\rightarrow M$ is the trajectory satisfying
the above property. Assume also that map $f=f_{d}$ is continuous. Then for
each $d\geq0$ we have a discrete time control system $D_{f}$.

Next, the system $C_{h}$ is called \textit{controllable} if any state can be
reached from any other state, i.e., for each pair of states $x,y\in M$ there
is a trajectory $c:[a,b]\rightarrow M$ such that $x=c(a),y=c(b).$

We make the same assumption about $f$ as in Section \ref{Controllability}: if
the initial state lies at the boundary $\partial M$ of $M$ then the next
state, regardless of the input, lies within a certain neighborhood $W$ of
$\partial M,$ or, alternatively, $f(\partial M\times U)\subset\partial M$. In
particular, this condition is satisfied if $h(x,u)$ is tangent to $\partial M$
for all $(x,u)\in Q$. Let $U^{\prime}$ be the set of controls that take any
given state to the the boundary $\partial M,$ i.e.,
\[
U^{\prime}=\{u\in U:f(x,u)\in\partial M\text{ for all }x\in M\}.
\]
Then $f$ is a map of pairs, $f:(M,\partial M)\times(U,U^{\prime}%
)\rightarrow(M,\partial M).$ Given a subset $L$ of $M,$ let $L^{\prime}%
=L\cap\partial M$ and $f^{\prime}:(L,L^{\prime})\times(U,U^{\prime
})\rightarrow(M,\partial M)$ be the restriction of $f.$

\begin{theorem}
\textbf{(Sufficient condition of controllability) }Suppose that there are
$a_{0}\in H_{p}(L,L^{\prime}),$ $v_{0}\in H_{s_{0}}(U,U^{\prime}),...,v_{r}\in
H_{s_{r}}(U,U^{\prime})$ such that
\[
a_{1}=f_{\ast}^{\prime}(a_{0}\otimes v_{0}),a_{2}=f_{\ast}(a_{1}\otimes
v_{1}),...,a_{r+1}=f_{\ast}(a_{r}\otimes v_{r})\neq0.
\]
Then the continuous time system $C_{h}$ is controllable from $L$ by means of
piece-wise constant controls.
\end{theorem}

\begin{proof}
The discrete system $D_{f}$ is controllable by Theorem \ref{SuffControl}.
\end{proof}

It follows also that given $\varepsilon>0$ if a map $k:Q\rightarrow TM$
satisfies $\rho(k(z),h(z))<\varepsilon$ for all $z\in Q$, where $\rho$ is the
distance on $TM,$ and the system $C_{k}$ satisfies all the above assumptions,
then $C_{k}$ is also controllable. We can say then that $C_{h}$ is
\textit{robustly controllable}$.$

Consider the applicability of this theorem to local controllability or
controllability in a Euclidean space. In either case $M$ is the $n$-ball. Its
homology (relative to the boundary) is nontrivial only in dimension $n$. As a
result the above "chain" of homology classes $a_{1},a_{2},...,a_{r+1}$ will
have to have only one "link", $a_{1}=f_{\ast}^{\prime}(a_{0}\otimes v_{0})\in
H_{n}(M,\partial M)\backslash\{0\}$. Thus the theorem reduces to the one-step
controllability provided $f_{\ast}^{\prime}:H_{n}(\{p\}\times(U,U^{\prime
}))\rightarrow H_{n}(M,\partial M)=\mathbf{Q}$ is nonzero. In particular, the
dimension of the space of inputs $U$ must be at least $n.$

Observe also that if $\partial M=\varnothing$, this theorem is vacuous.
Indeed, $f=f_{d}$ is homotopic to the constant map $f_{0}$ under the homotopy
$H(t,x,u)=f_{t}(x,u),$ hence $f_{\ast}=0$. Therefore the condition of the
theorem is never satisfied.

Here's another approach to controllability. Let $\mathcal{A}^{\prime}$ be the
set of controls whose trajectories have one of the end points at the boundary
of $M,$ i.e.,%
\[
\mathcal{A}^{\prime}=\{z:[0,d]\rightarrow Q,z\in\mathcal{A},c_{z}%
(0)\in\partial M\text{ or }c_{z}(d)\in\partial M\}.
\]
Define $G(u)=(c_{z}(0),c_{z}(d))$, the end points of the trajectory
$c_{z}=pz:[0,d]\rightarrow M$ corresponding to $z$. Then $G:(\mathcal{A}%
,\mathcal{A}^{\prime})\rightarrow(M\times M,\partial(M\times M))$ is a well
defined map of pairs.

\begin{theorem}
\textbf{\label{SuffControlCont}(Sufficient condition of controllability) }If
\[
G_{\ast}:H_{2n}(\mathcal{A},\mathcal{A}^{\prime})\rightarrow H_{2n}(M\times
M,\partial(M\times M))=\mathbf{Q}\text{ is non-zero}%
\]
then the continuous time system $C_{h}$ is controllable.
\end{theorem}

\begin{proof}
By Corollary \ref{Surj} $G$ is onto.
\end{proof}

A\ similar condition is found in \cite{Nistri}, where a boundary operator
$l:AC([0,1],\mathbf{R}^{n})\times L^{\infty}([0,1],\mathbf{R}^{n}%
)\rightarrow\mathbf{R}^{p}$ is considered instead of $G.$ One of the
conditions of controllability is $\deg l_{0}\neq0,$ where $l_{0}$ is the
restriction of $l$ to some $p$-dimensional subspace and $\deg l_{0}$ its
topological degree.

\href{http://inperc.com/wiki/index.php?title=Applications_of_Lefschetz_numbers_in_control_theory_by_Saveliev}{Lefschetz numbers in control theory}
\end{document}